\documentclass[fontsize=10pt]{scrreprt}

\usepackage[square,numbers,nonamebreak]{natbib}
\usepackage{setspace}

\usepackage[english]{babel}

\usepackage{chngcntr}
\counterwithout{equation}{section}

\usepackage[hang]{footmisc}
 at 11truept
\font\ssc=pplrc9d at 11 truept
\newcommand\qedbox{$\rlap{$\sqcap$}\sqcup$}

\usepackage[hmarginratio=1:1, bottom=2cm, top=2.5cm]{geometry}

\usepackage[automark]{scrlayer-scrpage}
\cohead{\textnormal{The structure skew brace is finitely presented}}
\lehead{\pagemark}
\rohead{\pagemark}

\let\ceheadL\cehead
\renewcommand\cehead[1]{
\ceheadL{\textnormal{#1}}
}

\cehead{Marco Trombetti}
\lefoot{}
\rofoot{}

\usepackage{graphicx}
\usepackage[usenames,dvipsnames,svgnames,table]{xcolor}
\definecolor{Maroon}{cmyk}{0, 0.87, 0.68, 0.32}
\definecolor{RoyalBlue2}{cmyk}{80,100,0,0.1}

\newcommand\auths[1]{\large \textsc{\textcolor{Maroon}{#1}}\setstretch{1.2}}
\newcommand\titl[1]{\center \linespread{1.1}\color{RoyalBlue2}\Large\textbf{ #1}\color{black}\bigskip} 
\renewcommand\abstract[1]{
\begin{center}
{\textbf{Abstract}}
\end{center}
{
\linespread{1.1}\fontsize{9pt}{-10pt}\selectfont #1}}

\usepackage[utf8]{inputenc}
\usepackage[sc,osf]{mathpazo}
\usepackage[euler-digits]{eulervm}
\usepackage[T1]{fontenc}
\usepackage{mathtools}

\DeclareSymbolFont{operators}{\encodingdefault}{ppl}{m}{n}
\DeclareMathAlphabet{\mathbf}{\encodingdefault}{ppl}{bx}{n}
\DeclareMathAlphabet{\mathit}{\encodingdefault}{ppl}{m}{it}

\usepackage{amsfonts}
\usepackage{amsmath}
\usepackage{ragged2e}
\usepackage{titlesec}
\usepackage{amssymb}

\renewcommand{\thesection}{\arabic{section}}
 
\titleformat{\section}{\medskip\bigskip\normalfont\Large\bf}{\thesection}{0.5em}{}
\titleformat{\subsection}{\smallskip\bigskip\normalfont\large\bf}{\thesubsection}{0.5em}{}

\usepackage{amsthm}
\newtheoremstyle{dotless}{}{}{\itshape}{}{\bfseries}{}{1em}{}

\theoremstyle{dotless}
\newtheorem{theo}{Theorem}

\newtheorem{lem}[theo]{Lemma}
\newtheorem{cor}[theo]{Corollary}

\renewenvironment{proof}{\smallbreak\noindent {\sc Proof \;---\;}}{\hfill\qedbox}

\counterwithout{theo}{section}
\numberwithin{theo}{section}

\DeclareOldFontCommand{\rm}{\normalfont\rmfamily}{\mathrm}
\DeclareOldFontCommand{\sf}{\normalfont\sffamily}{\mathsf}
\DeclareOldFontCommand{\tt}{\normalfont\ttfamily}{\mathtt}
\DeclareOldFontCommand{\bf}{\normalfont\bfseries}{\mathbf}
\DeclareOldFontCommand{\it}{\normalfont\itshape}{\mathit}
\DeclareOldFontCommand{\sl}{\normalfont\slshape}{\@nomath\sl}
\DeclareOldFontCommand{\sc}{\normalfont\scshape}{\@nomath\sc}

\usepackage{tikz}

\DeclareMathOperator{\Aut}{Aut}
\DeclareMathOperator{\op}{op}
\DeclareMathOperator{\Soc}{Soc}

\usepackage{mathtools}

\begin{document}

\setheadsepline{1pt}[\color{black}]

\titl{The structure skew brace associated with a finite non-degenerate solution of the Yang-Baxter equation is\\ finitely presented\footnote{The author is supported by GNSAGA (INdAM) and is a member of the non-profit association ``Advances in Group Theory and Applications'' (www.advgrouptheory.com).}}

\auths{Marco Trombetti}

\thispagestyle{empty}
\justify\noindent
\setstretch{0.3}
\abstract{The aim of this paper is to show that the structure skew brace associated with a finite non-degenerate solution of the Yang-Baxter equation is finitely presented.}

\setstretch{2.1}
\noindent
{\fontsize{10pt}{-10pt}\selectfont {\it Mathematics Subject Classification \textnormal(2020\textnormal)}: 16T25}\\[-0.8cm]

\noindent 
\fontsize{10pt}{-10pt}\selectfont  {\it Keywords}: Yang–Baxter equation, skew brace, free skew brace\\[-0.8cm]

\setstretch{1.1}
\fontsize{11pt}{12pt}\selectfont
\section{Introduction}

\noindent The Yang-Baxter equation (YBE) is one of the fundamental equations of physics. It takes its name from the independent works of the physicists Chen-Ning Yang \cite{Yang} and Rodney Baxter \cite{Baxter}. It turns out that this equation plays a relevant role in many different subjects such as knot theory, braid theory, operator theory, Hopf algebras, quantum groups, $3$-manifolds and the monodromy of differential equations. 

A {\it solution} to the YBE is a pair $(V,R)$, where $V$ is a vector space and $R$ is a linear map $R:\, V\otimes V\longrightarrow V\otimes V$ such that $$(R\otimes\operatorname{id})(\operatorname{id}\otimes R)(R\otimes\operatorname{id})=(\operatorname{id}\otimes R)(R\otimes\operatorname{id})(\operatorname{id}\otimes R).$$ At the present time, we are far from being able to provide a full classification of the solutions to the YBE. However, in recent years, there has been an assault at the so-called {\it set-theoretic} (or, {\it combinatorial}) solutions of the YBE, \hbox{i.e.} those solutions $(V,R)$ such that $R$ is induced by linear extensions of a bijective map $$r:\, X\times X\longrightarrow X\times X,$$ where $X$ is a basis of $V$ (see \cite{Drinfeld}); in this case, also the pair $(X,r)$ is called a {\it set-theoretic} (or, {\it combinatorial}) solution. This is mainly because the construction of set-theoretic solutions can sometimes be based on the use of (associative and non-associative) algebraic structures. Among these structures, skew braces hold a prominent position. The theory of skew braces arises from the theory of Jacobson radical rings and can be found in many contexts in mathematics (see for example \cite{Vendramin}). In \cite{Rump2}, Wolfgang~Rump found out that radical rings give rise to {\it involutive non-degenerate} solutions, i.e. solutions $(X,r)$ such that $r^2=\operatorname{id}$, and if we write $$r(x,y)=(\sigma_x(y),\tau_y(x))$$ then the maps $\sigma_x$ and $\tau_y$ are bijective, for all $x,y\in X$. There, radical rings have been generalised to braces, and it was proved that all involutive non-degenerate solutions of the YBE are restrictions of those solutions obtained from left braces; actually, radical rings correspond to two-sided braces (see \cite{Rump2}), \hbox{i.e.} left braces that are also right braces. A brace is essentially a set endowed with two group structures one of which is abelian, linked together with a ``distributivity-like'' relation (see next section for the precise definitions). By removing the abelianity constraint, one obtains the concept of a skew brace (see \cite{OriginalBrace}), which is a powerful tool for studying  (not necessarily involutive) non-degenerate solutions. 

If $(X,r)$ is a non-degenerate solution, the skew brace controlling the structure of this solution is the {\it structure skew brace} $G(X,r)$. Understanding the structure of $G(X,r)$ makes it possible to better understand the structure of the corresponding solution.

Free skew braces on a given set $X$ exist by a general result in universal algebra (see for instance \cite{bergman}). A simpler construction has recently been given in \cite{free}. It therefore makes sense to speak of a presentation of a skew brace. The aim of this short paper is to prove the following theorem.

\medskip

\noindent{\bf Main Theorem}\quad {\it Let $(X,r)$ be a finite non-degenerate solution of the YBE. Then the structure skew brace $B=G(X,r)$ is finitely presented.}

\medskip

The main step in the proof is Theorem \ref{thfp}, which is of an independent interest (see, for example, Theorem \ref{alternative}) and shows how to extend finite presentations of skew braces.

\section{Preliminaries and notation}

\noindent Let $B$ be a set. If $(B,+)$ and $(B,\circ)$ are (not necessarily abelian) groups, then the triple $(B,+,\circ)$ is a {\it skew \textnormal(left\textnormal) brace} if the skew (left) distributive property $$a\circ(b+c)=a\circ b-a+a\circ c$$ holds for all $a,b,c\in B$. Now, let $(B,+,\circ)$ be a skew brace. We refer to $(B,+)$ as the \emph{additive group} of $B$ and to $(B,\circ)$ as the \emph{multiplicative group} of $B$. We denote by $0$ the identity of $(B,+)$, by $1$ the identity of $(B,\circ)$, and  by $-a$ and $a^{-1}$ the inverses of $a$ in $(B,+)$ and $(B,\circ)$, respectively. The skew distributive property implies $0 = (-1)\circ 0 = (-1)\circ (0+0) = (-1) \circ 0 -(-1) + (-1)\circ 0 = 1$, i.e. $0=1$. It should be also noticed that the map $$\lambda:\, a\in (B,\circ)\mapsto \left(\lambda_a:\, b\mapsto \lambda_a(b)=-a+a\circ b\right)\in\Aut(B,+)$$ is a group homomorphism and the following relations hold $$a+b=a\circ\lambda_a^{-1}(b),\quad a\circ b=a+\lambda_a(b),\quad -a=\lambda_a\left(a^{-1}\right).$$ 

\indent In analogy with ring theory, a third relevant (non-necessarily associative) operation in skew braces is defined as follows $$a\ast b=\lambda_a(b)-b=-a+a\circ b-b$$ and one easily checks that it satisfies the relations 
\[
\begin{array}{c}
a \ast (b + c) = a \ast b + b + a \ast c - b\quad\textnormal{and}\\[0.2cm]
(a \circ b) \ast c=a \ast (b \ast c) + b \ast c + a \ast c
\end{array}
\] 
for all $a,b,c\in B$. If the additive group $(B,+)$ of $B$ is abelian, we call $B$ simply a {\it \textnormal(left\textnormal) brace}, or, a {\it \textnormal{(}left\textnormal{)} brace of abelian type}. If $(G,\cdot)$ is any group, then $(G,\cdot,\cdot)$ is a skew brace called a {\it trivial skew brace} and $(G,\cdot^{\op},\cdot)$ is a skew brace called an \emph{almost trivial skew brace}; if $(G,\cdot)$ is abelian, then the trivial skew brace and the almost trivial skew brace coincide and we simply speak of a {\it trivial brace}.

If we consider the natural semidirect product $G=(B,+)\rtimes(B,\circ)$, where \[
\begin{array}{c}
(a,b)(c,d)=(a + \lambda_b(c),b\circ d)
\end{array}
\] for all $a,b,c,d\in B$, then an easy computation shows that the $\ast$-operation corresponds to a commutator of type \[
\begin{array}{c}\label{equationfinale}
\left[(0,a),(b,0)\right]=(a\ast b,0),
\end{array}
\] for all $a,b\in B$ (notice that our convention for commutators in a group $(G,\cdot)$ is the following one: $[x,y]=xyx^{-1}y^{-1}$).

A {\it left ideal} of a skew brace $B$ is a subgroup $I$ of $(B,+)$ such that $\lambda_a(I)\subseteq I$ for all $a\in B$; this is equivalent to $B\ast I\subseteq I$, so $I$ is also a subgroup of $(B,\circ)$. An {\it ideal} is a left ideal that is normal in $(B,+)$ and $(B,\circ)$ (notice that the last condition is equivalent to demanding that $I\ast B\subseteq I$); in this case, it is known that $B/I$ is a skew brace and $a+I=a\circ I$ for all $a\in B$. A skew brace is {\it simple} if it has no proper non-zero ideals.

The {\it socle} of $B$ is defined as $\Soc(B)=\operatorname{Ker}(\lambda)\cap Z(B,+)$ and the {\it annihilator} (see \cite{catino}) of $B$ is defined as $\operatorname{Ann}(B)=\operatorname{Soc}(B)\cap Z(B,\circ)$. Moreover, we let $B^{(2)}=B\ast B$ be the subgroup of $(B,+)$ generated by all elements of the form $a\ast b$ for all $a,b\in B$. It can be proved that $\operatorname{Soc}(B)$, $\operatorname{Ann}(B)$ and $B^{(2)}$ are ideals. In connection with $B^{(2)}$, we further observe that $B$ is a trivial skew brace if and only if $B^{(2)}=\{0\}$.

A {\it system of generators} for the skew brace $B$ is just a subset $S$ of $B$ such that $B$ is the smallest sub-skew brace of $B$ containing $S$. As usual, if $B$ has a finite system of generators, we say that $B$ is {\it finitely generated}. 

The following definitions have been introduced in \cite{FCbraces} in order to deal with finitely generated skew braces. Let $\mathcal{S}=\{x_1,\ldots,x_n\}$ be symbols. A {\it $b$-word} with respect to~$\mathcal{S}$ is a sequence of symbols recursively defined as follows: the empty sequence is a $b$-word and such are the $1$-element sequences $x_1$, $x_2$, \ldots, $x_n$; if we have two $b$-words~$w_1$ and~$w_2$, then the sequences $w_1\circ w_2$, $w_1+w_2$, $-w_1$, $w_1^{-1}$ are $b$-words.  Now, let~$B$ be a skew brace and let $b_1,\ldots, b_n$ be elements of $B$. It is clear that if $w(x_1,\ldots,x_n)$ is any $b$-word, then we may evaluate $w(b_1,\ldots,b_n)$ in $B$. Thus, the smallest sub-skew brace~$C$ generated by $b_1,\ldots, b_n$ in $B$ is precisely the set of all evaluations of $b$-words with respect to $b_1,\ldots,b_n$. Note also that we usually abbreviate an evaluation like $w(b_1,\ldots,b_n)$ by simply writing $w(b)$.

Finally, we illustrate the connection between skew braces and solutions of the YBE. Let $B$ be a skew brace and let $$r_B:\, (a,b)\in B\times B\mapsto \left(\lambda_a(b),\,\lambda_a(b)^{-1}\circ a\circ b\right)\in B\times B.$$ Then $(B,r_B)$ is a non-degenerate solution of the YBE. Conversely, if $(X,r)$ is a non-degenerate solution of the YBE, then there is a unique skew brace structure over the structure group $$G(X,r)=\langle X\, |\, xy=\sigma_x(y)\tau_y(x),\; x,y\in X\rangle$$ such that $r_{G(X,r)}(\iota\times\iota)=(\iota\times\iota)r$, where $\iota:\, X\longrightarrow G(X,r)$ is the canonical map: the multiplicative group of this skew brace is $G(X,r)$ and the additive is $$A(x,r)=\langle X\,|\, x+\sigma_x(y)=\sigma_x(y)+\sigma_{\sigma_x(y)}\big(\tau_y(x)\big)\;\textnormal{for all}\; x,y\in X\rangle$$ (see \cite{Soloviev},\cite{Zhu}). We refer to this skew brace as the {\it structure skew brace} of $(X,r)$.

The second skew brace associated with a non-degenerate solution of the YBE is the permutation skew brace $\mathcal{G}(X,r)$, that is, the subgroup of $\operatorname{Sym}(X)\times\operatorname{Sym}(X)$ generated by the elements of the form $(\sigma_x,\tau_x^{-1})$, with $x\in X$. It turns out that $\mathcal{G}(X,r)$ is a homomorphic image of the structure skew brace $G(X,r)$ and the kernel of this homomorphism is contained in $\operatorname{Soc}\big(G(X,r)\big)$. In particular, if $X$ is finite, then $G(X,r)/\operatorname{Soc}\big(G(X,r)\big)$ is finite and the additive group of $G(X,r)$ is central-by-finite.

\section{Proof of the Main Theorem}


Let $B$ be a skew brace. A {\it presentation} of $B$ is an exact sequence of skew braces \[
\begin{array}{c}\label{exact}
0\rightarrow R\rightarrow F\xrightarrow{\theta} B\rightarrow 0,\tag{$\star$}
\end{array}
\] where $F$ is a free brace (see for instance \cite{free}) over some set $X$. Suppose $B$ is finitely generated as a skew brace, so $X$ may be assumed to be finite and we may write $X=\{x_1,\ldots,x_m\}$; clearly $B$ is generated by $a_1=\theta(x_1),\ldots,a_m=\theta(x_m)$.

The exact sequence \eqref{exact} is a {\it finite presentation} of $B$ if there are finitely many elements $\rho_1,\ldots,\rho_n$ of $F$ such that $\operatorname{Ker}(\theta)$ is the ideal of $F$ generated by $\rho_1,\ldots,\rho_n$. In this case, the skew brace $B$ is said to be {\it finitely presented} by $a_1,\ldots,a_m$ subject to the relations $\rho_1=\ldots=\rho_n=1$.

We show that this definition is independent of the presentation in the sense that if $b_1,\ldots,b_l$ is another finite set of generators for $B$, then $B$ can be presented by the $b_i$ subject to a set of $l+n$ relations. For suitable $b$-words $\theta_i$ and $\phi_j$ we can write $$a_i=\theta_i(b)\quad\textnormal{and}\quad b_j=\phi_j(a).$$ Then the relations \begin{equation}\label{seconda}
\rho_k\big(\theta_1(b),\ldots,\theta_m(b)\big)=1\quad\textnormal{and}\quad b_j=\phi_j\big(\theta_1(b),\ldots,\theta_m(b)\big),
\end{equation}
 where $k=1,\ldots,n$ and $j=1,\ldots,l$, are satisfied in $B$. Let $\overline{B}$ be the skew brace presented by a set $\{\overline{b}_1,\ldots,\overline{b}_l\}$ subject to the $l+n$ relations \eqref{seconda}. Since all the defining relations of $\overline{B}$ are satisfied by the $b_j$, the map $\overline{b}_j\mapsto b_j$ defines a homomorphism $\beta$ of $\overline{B}$ onto $B$. Let $\overline{a}_i=\theta_i(\overline{b})$, so $\overline{B}$ can also be generated by $\overline{a}_1,\ldots,\overline{a}_m$ since $\overline{b}_j=\phi_j(\overline{a}_1,\ldots,\overline{a}_m)$. Since the original defining relations for $B$ in the $a_i$ are also valid in the $\overline{a}_i$, the map $a_i\mapsto\overline{a}_i$ determines a homomorphism $\alpha$ of $B$ onto $\overline{B}$. Now, $\alpha\beta$ and $\beta\alpha$ are identity maps, so~$\alpha$ and $\beta$ are isomorphisms and $B\simeq\overline{B}$. Hence $B$ is presented by the $b_j$ subject to the relations \eqref{seconda}.

\begin{lem}\label{fgideal}
Let $I$ be an ideal of a skew brace $B$. If $B$ is finitely generated and $B/I$ is finitely presented, then $I$ is finitely generated as an ideal of $B$.
\end{lem}
\begin{proof}
Let $B$ be generated by $b_1,\ldots,b_m$ and let $F$ be the free skew brace on $\{x_1,\ldots,x_m\}$. Let $\theta:\, F\longrightarrow B$ be the homomorphism mapping $x_i$ to $b_i$, for each $i$. Let $J=\theta^{-1}(I)$. Then the mapping $x_i\mapsto b_i+I$ defines a homomorphism of $F$ onto $B/I$ with kernel $J$. Now, the proof before the statement shows that there are elements $\rho_1,\ldots,\rho_n$ of $F$ such that $J$ is the ideal of $F$ generated by $\rho_1,\ldots,\rho_n$. Therefore $I=\theta(J)$ is generated as an ideal of $B$ by $\theta(\rho_1),\ldots,\theta(\rho_n)$. The statements is proved.
\end{proof}

\medskip

Among finitely presented skew braces there certainly are finite skew braces (just put as relations those deduced from its multiplicative and additive tables). Other classes of finitely presented skew braces can be obtained through the following result.

\begin{theo}\label{thfp}
Let $I$ be an ideal of a skew brace $B$ such that \begin{itemize}
\item[\textnormal{(1)}] $(I,+)$, $(I,\circ)$, $(B/I,\circ)$ are finitely generated;
\item[\textnormal{(2)}] $I$, $B/I$, $(B/I,+)$ are finitely presented (the last one as a group).
\end{itemize}

\noindent Then $B$ is finitely presented.
\end{theo}
\begin{proof}
Let $b_1,\ldots,b_m$ be elements of $B$ such that $(B/I,+)$ and $(B/I,\circ)$  are both generated as semigroups by $b_1+I=b_1\circ I$, \ldots, $b_m+I=b_m\circ I$: it is enough to put together the representatives of the generators of $(B/I,+)$ and their additive inverses, and the representatives of the generators of $(B/I,\circ)$ and their multiplicative inverses. Similarly, let $0=a_1,\ldots,a_n$ be elements of $I$ generating both $(I,+)$ and $(I,\circ)$ as sub-semigroups.

Let $\rho_1=\ldots=\rho_\ell=1$ be the relations to which $a_1,\ldots,a_n$ are subject in generating~$I$ as a skew brace. For each $i,j=1,\ldots,m$, write \[
\begin{array}{c}
  b_i\circ b_j=\sum_{k=1}^{\sigma_\circ(i,j)}b_{(k;i,j,\circ)}+\mu_{i,j}^\circ(a_1,\ldots,a_n)\quad\textnormal{and}\\[0.5cm]
 b_i^{-1}=\prod_{k=1}^{\sigma_1(i,j)}b_{(k;i,1)}+\mu_{i}^1(a_1,\ldots,a_n),\;
\end{array}\]   for certain $b$-words $\mu_{i,j}^\circ$, $\mu_i^1$ and $(k;i,j,\circ),(k;i,1)\in\{1,\ldots,m\}$. Moreover, for all $i=1,\ldots,m$ and $j=1,\ldots,n$, write \[
\begin{array}{c}
\lambda_{b_i}(a_j)=\theta_{i,j}^1(a_1,\ldots,a_n),\\[0.3cm]
a_j^{\circ, b_i}=\theta_{i,j}^2(a_1,\ldots,a_n),\\[0.3cm]
a_j^{+, y_i}=\theta_{i,j}^3(a_1,\ldots,a_n)\;\,\textnormal{and}\\[0.3cm]
\lambda_{b_i^{-1}}(a_j)=\theta_{i,j}^4(a_1,\ldots,a_n),\\[0.3cm]
\end{array}
\] for certain $b$-words $\theta_{i,j}^1$, $\theta_{i,j}^2$, $\theta_{i,j}^3$, $\theta_{i,j}^4$. Let $\sigma_1=\ldots=\sigma_s=1$ be the relations to which\linebreak $b_1+I,\ldots,b_m+I$ are subject in generating $B/I$ as a skew brace; write \[
\begin{array}{c}
\sigma_1(b_1,\ldots,b_m)=\eta_1^1(a_1,\ldots,a_n),\;\ldots,\;\sigma_m(b_1,\ldots,b_m)=\eta_m^1(a_1,\ldots,a_n),
\end{array}
\] for certain $b$-words $\eta_i$. Let $\tau_1=\ldots=\tau_r=1$ be the relations to which \hbox{$b_1+I,\ldots,b_m+I$} are subject in generating $(B/I,+)$; write \[
\begin{array}{c}
\tau_1(b_1,\ldots,b_m)=\eta_1^2(a_1,\ldots,a_n),\;\ldots,\;\tau_r(b_1,\ldots,b_m)=\eta_r^2(a_1,\ldots,a_n),\\
\end{array}
\] for certain words $\eta_i^2$.

Let $F$ be the free skew brace on $\{y_1,\ldots,y_m,x_1,\ldots,x_n\}$. Let $J$ be the ideal of $F$ generated by the following elements \begin{equation}\label{unsacco}
\begin{cases}
y_i\circ y_j-\sum_{k=1}^{\sigma_\circ(i,j)}y_{(k;i,j,\circ)}-\mu_{i,j}^\circ(x_1,\ldots,x_n) & (i,j=1,\ldots,m)\\
y_i^{-1}-\prod_{k=1}^{\sigma_1(i,j)}y_{(k;i,1)}-\mu_i^1(x_1,\ldots,x_n) & (i=1,\ldots,m)\\
\rho_i(x_1,\ldots,x_n) & (i=1,\ldots,\ell)\\
\lambda_{y_j}(x_i)-\theta_{i,j}^1(x_1,\ldots,x_n) & (i=1,\ldots,n;\, j=1,\ldots,m)\\
x_i^{\circ,y_j}-\theta_{i,j}^2(x_1,\ldots,x_n) & (i=1,\ldots,n;\, j=1,\ldots,m)\\
x_i^{+,y_j}-\theta_{i,j}^3(x_1,\ldots,x_n) & (i=1,\ldots,n;\, j=1,\ldots,m)\\
\lambda_{y_j^{-1}}(x_i)-\theta_{i,j}^4(x_1,\ldots,x_n) & (i=1,\ldots,n;\, j=1,\ldots,m)\\
\sigma_i(y_1,\ldots,y_m)-\eta_i^1(x_1,\ldots,x_n) & (i=1,\ldots,s)\\
\tau_i(y_1,\ldots,y_m)-\eta_i^2(x_1,\ldots,x_n) & (i=1,\ldots,r)\\
\end{cases}
\end{equation}
Let $\overline{F}=F/J$. The assignation $x_i\mapsto a_i$ and $y_j\mapsto b_j$ defines a homomorphism $\varphi$ of $\overline{F}$ over $B$. Let $\overline{C}$ be the sub-skew brace of $\overline{F}$ generated  by $\overline{x}_1,\ldots,\overline{x}_n$. Then the restriction of~$\varphi$ to $\overline{C}$ is injective by the third of \eqref{unsacco}, so $\overline{C}$ and $I$ are isomorphic skew braces and in particular $\overline{C}$ is generated both additively and multiplicatively by $x_1,\ldots,x_n$.

Let $\overline{c}\in \overline{C}$ and $j\in\{1,\ldots,m\}$. Let $\overline{y}_1',\ldots,\overline{y}_h'$ in $\{\overline{y}_1,\ldots,\overline{y}_m\}$. Then \begin{equation}\label{continuo}
\begin{array}{c}
\overline{y}_1'\circ\ldots\circ\overline{y}_h'\circ\overline{c}=\overline{y}_1'\circ\ldots\circ\overline{y}_h'\circ(0+\overline{c})=\overline{y}_1'\circ\ldots\circ\overline{y}_{h-1}'\circ(\overline{y}_h'+\lambda_{\overline{y}_h'}(\overline{c}))\\[0.2cm]
=\ldots=\overline{y}_1'\circ\ldots\circ\overline{y}_h'+\lambda_{\overline{y}_1'}\big(\lambda_{\overline{y}_2'}\big(\ldots\big(\lambda_{\overline{y}_h'}(\overline{c})\big)\big)\ldots\big)
\end{array}
\end{equation} Write $\overline{c}=\overline{x}_1'+\ldots+\overline{x}_r'$ for certain $\overline{x}_i'\in\{\overline{x}_1,\ldots,\overline{x}_n\}$. By the fourth of \eqref{unsacco}, we have that $\lambda_{\overline{y}_1'}\big(\lambda_{\overline{y}_2'}\big(\ldots\big(\lambda_{\overline{y}_h'}(\overline{c})\big)\big)\ldots\big)$ belongs to $\overline{C}$. Moreover, the seventh of \eqref{unsacco} shows that the equation  $$\lambda_{\overline{y}_1'}\big(\lambda_{\overline{y}_2'}\big(\ldots\big(\lambda_{\overline{y}_h'}(\overline{c})\big)\big)\ldots\big)=\overline{d},$$ where $\overline{d}$ is in $\overline{C}$, has always a solution in $\overline{C}$. This means that every element of the form $\overline{y}_1'\circ\ldots\circ\overline{y}_h'+\overline{d}$, where $\overline{d}\in\overline{C}$, can always be written as an element of type $\overline{y}_1'\circ\ldots\circ\overline{y}_h'\circ\overline{c}$ for some $\overline{c}\in\overline{C}$, and vice-versa.

Now, using the first of \eqref{unsacco}, the final side of \eqref{continuo} can be written as $$\overline{y}_1''+\ldots+\overline{y}_k''+\overline{d}$$ for a certain $\overline{d}\in\overline{C}$ and elements $\overline{y}_i''\in\{\overline{y}_1,\ldots,\overline{y}_m\}$. Conversely, assume we have an element of the form $$\overline{b}=\overline{y}_1''+\ldots+\overline{y}_k''+\overline{d}$$ for some $\overline{d}\in\overline{C}$. Let $$Y_+=\langle\overline{y}_1,\ldots,\overline{y}_m,\overline{x}_1,\ldots,\overline{x}_n\rangle_+$$ be the additive group  generated by $\overline{y}_1,\ldots,\overline{y}_m,\overline{x}_1,\ldots,\overline{x}_n$. It follows from the sixth of \eqref{unsacco} that $\overline{C}$ is normal in $Y$. Consequently, $Y_+/\overline{C}$ is isomorphic to $(B/I,+)$ by the ninth of \eqref{unsacco}, and hence $\varphi$ induces an isomorphism of $Y_+$ onto $(B,+)$. Of course, $\varphi(\overline{b})$ can be written as a product of elements in $\{b_1,\ldots,{b}_m,{a}_1,\ldots,{a}_n\}$; and the corresponding product $\overline{b}'$, in terms of $\{\overline{y}_1,\ldots,\overline{y}_m,\overline{x}_1,\ldots,\overline{x}_n\}$, belongs to $Y_+$ by what we have already proved. This means that $\overline{b}=\overline{b}'$, so $\overline{y}_1''+\ldots+\overline{y}_k''+\overline{d}$ can be written as a product \hbox{$\overline{y}\circ\overline{e}$,} where $\overline{y}$ is a product of elements in $\{\overline{y}_1,\ldots,\overline{y}_m\}$ and $\overline{e}\in\overline{C}$, and vice-versa.

As a by-product of the previous argument, we also find that, for every $i=1,\ldots,m$, a relation of type $$-\overline{y}_i+\overline{C}=\sum_{k=1}^{\sigma_2(i,j)}\overline{y}_{(k;i,2)}+\overline{C},$$ where $(k;i,2)\in\{1,\ldots,m\}$, holds.

Now, we claim that every element of $\overline{F}$ can be written as a sum $\overline{y}+\overline{c}$, where $\overline{y}$ is a sum of elements in $\{\overline{y}_1,\ldots,\overline{y}_m\}$  and $\overline{c}\in\overline{C}$. Let $$\overline{E}=\left\{\sum_{j=1}^{k}\overline{y}_{i_j}+\overline{c}\,:\,\overline{c}\in\overline{C}\quad\textnormal{and}\quad i_j\in\{1,\ldots,m\}\right\}.$$ Of course, using what we have proved, we get that $$\overline{E}=\left\{\prod_{j=1}^{k}\overline{y}_{i_j}\circ\overline{c}\,:\,\overline{c}\in\overline{C}\quad\textnormal{and}\quad i_j\in\{1,\ldots,m\}\right\}.$$ Also using the remaining relations of \eqref{unsacco}, it is very easy to see that $\overline{E}$ is a sub-skew brace of $\overline{F}$, and consequently that $\overline{E}=\overline{F}$. The claim is proved.

The previous claim yields that $\overline{C}$ is an ideal of $\overline{F}$. Clearly, $\overline{F}/\overline{C}$ is isomorphic to~$B/I$ by the last of \eqref{unsacco} and consequently $\varphi$ is an isomorphism. Therefore $B$ is finitely presented.~\end{proof} 

\medskip

Before proving our main theorem, we deal with some of the most relevant consequences of the previous result. Recall first that a skew brace is said to satisfy {\it ACC} if it satisfies the ascending chain condition on sub-skew braces. If $I$ is an ideal of a skew brace $B$ such that both $I$ and $B/I$ satisfy ACC, then $B$ satisfy ACC by a standard argument. For trivial skew braces, ACC is just the usual ascending chain condition on subgroups.

\begin{cor}\label{maxbaer}
Let $B$ be a skew brace having an ascending series of ideals $$I_0=\{0\}\subseteq I_1\subseteq\ldots I_\alpha\subseteq I_{\alpha+1}\subseteq\ldots I_\lambda=B$$ \textnormal(here $\lambda$ is an ordinal number and $\alpha<\lambda$\textnormal) such that for each $\alpha<\lambda$:

\begin{itemize}
\item[\textnormal{(i)}] $(I_{\alpha+1}/I_\alpha,\circ)$ is finitely generated;
\item[\textnormal{(ii)}] $(I_{\alpha+1}/I_\alpha,+)$ and $I_{\alpha+1}/I_\alpha$ are finitely presented;
\item[\textnormal{(iii)}] $(I_{\alpha+1}/I_\alpha,+)$ and $I_{\alpha+1}/I_\alpha$ satisfy \textnormal{ACC}.
\end{itemize} 

\noindent Then every finitely generated sub-skew brace of $B$ satisfies \textnormal{ACC}.
\end{cor}
\begin{proof}
Let $\mathcal{I}$ be the set of all ideals $I$ of $B$ with the following property: 

\begin{itemize}
\item if $J$ is a sub-skew brace of $B$ containing $I$, and such that $J/I$ satisfies (i),(ii) and (iii), replacing $I_{\alpha+1}/I_\alpha$ by $J/I$, then every finitely generated sub-skew brace of $J$ satisfies ACC.
\end{itemize}

We claim that the set $\mathcal{I}$ is inductive. Let $J=\bigcup_{u\in U} J_u$ be the union of a chain $\{J_u\}_{u\in U}$ in $\mathcal{I}$. Assume that $K$ is a sub-skew brace of $B$ containing $J$ and such that $K/J$ satisfies (i),(ii) and (iii), mutatis mutandis. Let $C$ be a finitely generated sub-skew brace of $K$. Of course, we may assume $C+J=K$, because $K/J$ is finitely generated and we only need to show that $C$ satisfies ACC. Now, $$K/J=(C+J)/J\simeq C/(C\cap J)$$ is finitely presented, so Lemma \ref{fgideal} yields that $C\cap J$ is finitely generated as an ideal of $C$. But $C\cap J=\bigcup_{u\in U} (C\cap J_u)$ and so there is $\overline{u}\in U$ such that $C\cap J=C\cap J_{\overline{u}}$. This means that $$(C+J_{\overline{u}})/J_{\overline{u}}\simeq C/(C\cap J)$$ satisfies (i),(ii) and (iii), mutatis mutandis. Since $J_{\overline{u}}\in\mathcal{I}$, we have that $C$ satisfies ACC. This proves that $\mathcal{I}$ is inductive, so by Zorn's lemma it admits a maximal element, say~$M$. 

If $M=B$, we are done. Assume $M\neq B$. Since the class of finitely presented (trivial) skew braces satisfying ACC is closed with respect to forming quotients, our hypothesis implies that $B/M$ contains a non-zero ideal $N/M$ satisfying (i),(ii) and (iii), mutatis mutandis.  Then by Theorem \ref{thfp} (and the extension closure of ACC), $N$ belongs to $\mathcal{I}$, contradicting the maximality of $M$. The statement is proved.
\end{proof}

\medskip

The above result can be useful in the study of the structure of certain relevant types of skew braces. For example, Theorem \ref{alternative} provides an alternative proof of Corollary~3.4 of \cite{JeVAV22x} with an additional item.

\begin{cor}\label{trivialfp}
Let $B$ be a trivial brace such that $(B,+)$ is finitely generated. Then $B$ is finitely presented.
\end{cor}
\begin{proof}
Since $(B,+)$ is finitely generated, it can be written as the sum of $n$ cyclic groups $\langle b_1\rangle$, \ldots, $\langle b_n\rangle$. We prove the result by induction on $n$. 

Assume $n=1$. If $B$ is finite, the statement is clear. Thus we may also assume $B$ is infinite. Let $F$ be the free skew brace on $X=\{x_1\}$ and let $\varphi:\, F\longrightarrow B$ be the homomorphism mapping $x_1$ to $b_1$. Moreover, let $I$ be the ideal of $F$ generated by $x_1\ast x_1$, $x_1\ast (-x_1)$, $(-x_1)\ast(-x_1)$ and put $\overline{F}=F/I$. Then $\overline{x}_1+\overline{x}_1=\overline{x}_1\circ\overline{x}_1$, $$\overline{x}_1\circ(\overline{x}_1\circ\overline{x}_1)=\overline{x}_1\circ(\overline{x}_1+\overline{x}_1)=\overline{x}_1\circ\overline{x}_1-\overline{x}_1+\overline{x}_1\circ\overline{x}_1=\overline{x}_1+\overline{x}_1+\overline{x}_1,$$ and, more generally, $$\overline{x}_1^\ell=\ell \overline{x}_1$$ for every positive $\ell$. Since $0=\overline{x}_1-\overline{x}_1=\overline{x}_1\circ(-\overline{x}_1)$, we have that $-\overline{x}_1=\overline{x}_1^{-1}$ and consequently (as above) $$\overline{x}_1^\ell=\ell \overline{x}_1$$ for every (negative) $\ell$. Therefore $\overline{F}$ is isomorphic to $B$ and we are done.

Assume $n>1$. Let $I$ be the subgroup of $B$ generated by $b_1$. By induction, $B/I$ is finitely presented, and of course also $(B/I,+)$ is finitely presented. Thus Theorem \ref{thfp} yields that $B$ is finitely presented and completes the proof.
\end{proof}

\medskip

Let $B$ be a skew brace. We recursively define the {\it upper annihilator series} of~$B$ as follows (see \cite{FCbraces}). Put $\operatorname{Ann}_0(B)=\{0\}$ and $\operatorname{Ann}_1(B)=\operatorname{Ann}(B)$. If~$\alpha$ is an ordinal number, put $$\operatorname{Ann}_{\alpha+1}(B)/\operatorname{Ann}_\alpha(B)=\operatorname{Ann}\big(B/\operatorname{Ann}_{\alpha}(B)\big).$$ If $\mu$ is a limit ordinal, put $$\operatorname{Ann}_\mu(B)=\bigcup_{\gamma<\mu}\operatorname{Ann}_\gamma(B).$$ The last term of the upper socle series is the {\it hyper-annihilator} of $B$ and is denoted by $\overline{\operatorname{Ann}}(B)$. Note that every term of the upper annihilator series is an ideal of $B$ and it is easy to see that $\overline{\operatorname{Ann}}\big(B/\overline{\operatorname{Ann}}(B)\big)$ is trivial. If $B=\overline{\operatorname{Ann}}(B)$, then $B$ is called {\it annihilator hypercentral}; if $B=\operatorname{Ann}_n(B)$ for some positive integer $n$, then $B$ is called {\it annihilator nilpotent} (see \cite{Bonatto} and \cite{MR4256133}).

\begin{theo}\label{teoooo}
Let $B$ be a skew brace. If $B$ is annihilator hypercentral, then every finitely generated skew brace of $B$ satisfies ACC. 
\end{theo}
\begin{proof}
By Corollary \ref{trivialfp}, the upper annihilator series of $B$ can be refined to an ascending series satisfying the hypothesis  of Corollary \ref{maxbaer}. The statement is proved.~\end{proof}

\begin{cor}\label{corollario}
Let $B$ be a finitely generated skew brace. If $B$ is annihilator hypercentral, then~$B$ is annihilator nilpotent.
\end{cor}

\begin{theo}\label{alternative}
Let $B$ be an annihilator hypercentral skew brace. The following statements are equivalent:
\begin{itemize}
\item[\textnormal{(1)}] $B$ is finitely generated;
\item[\textnormal{(2)}] $B$ is finitely presented;
\item[\textnormal{(3)}] $(B,+)$ is finitely generated;
\item[\textnormal{(4)}] $(B,\circ)$ is finitely generated.
\end{itemize}
\end{theo}
\begin{proof}
It is clear that (3) implies (1), (4) implies (1) and (2) implies (1). Assume~(1) and let $n$ be such that $\operatorname{Ann}_n(B)=B$ (see Corollary \ref{corollario}). If $n\leq1$, the statement is just Corollary \ref{trivialfp}. Assume $n>1$. It follows from Lemma \ref{fgideal} that $\operatorname{Ann}(B)$ is finitely generated as an ideal of $B$, so it is finitely generated both as a multiplicative and an additive group, and it is also finitely presented by Corollary \ref{trivialfp}. Moreover, $(B/I,+)$ and $(B/I,\circ)$ are finitely generated nilpotent groups, so they are finitely presented. It follows that $(B,+)$ and $(B,\circ)$ are finitely generated.

Finally, since by induction $B/I$ is finitely presented, it follows from Theorem \ref{thfp} that~$B$ is finitely presented.
\end{proof}

\medskip

In Theorem \ref{alternative}, the hypothesis ‘‘annihilator hypercentral’’ can be replaced (essentially with the same proof) by weaker ones: for example, you may ask that $B$ has an ascending series of ideals $$\{0\}=I_0\subseteq\ldots I_\alpha\subseteq I_{\alpha+1}\subseteq\ldots I_\lambda=B$$ whose factors $I_{\alpha+1}/I_\alpha$ are either finite or contained in $\operatorname{Ann}(B/I_{\alpha})$.

\medskip

Our final result applies in particular to the structure skew brace associated with a finite non-degenerate solution of the YBE.  In fact, as we noted in the preliminaries, the socle of a such a skew brace has finite index; moreover, the multiplicative group is virtually abelian by  \cite[Corollary~6.2]{LeVe19}. This proves our main theorem. 

Theorem \ref{nilpotent2} also says something more on the structure skew brace in connection with a recently investigated class of skew braces.  Let $B$ be a skew brace. If $x\in B$, any element of type $g\ast x$, $x\ast g$, $g\circ x\circ g^{-1}$, $g+x-g$, for some $g\in B$, will be referred to as a {\it conjugate} of~$x$. Skew braces in which every element has finitely many conjugates have recently been investigated in \cite{FCbraces}; these skew braces were referred to as {\it skew braces satisfying property}~(S). In \cite{FCbraces} is proved that not every structure skew brace has property~(S) and that these skew braces are strictly connected with finite solution whose derived solution is indecomposable. As a by-product, the following result shows that every structure skew brace associated with a finite non-degenerate solution of the~YBE is actually very close to satisfy property (S).

\begin{theo}\label{nilpotent2}
Let $B$ be a finitely generated skew brace such that $B/\operatorname{Soc}(B)$ is finite. If $(B,\circ)$ is virtually abelian, then $B$ is finitely presented. Moreover, $B$ contains an ideal $I$ such that: 
\begin{itemize}
\item[\textnormal{(1)}] $I$ is a finitely generated trivial brace;
\item[\textnormal{(2)}] $I$ is contained in $\operatorname{Soc}(B)$;
\item[\textnormal{(3)}] $B/I$ is finite;
\item[\textnormal{(4)}] every element of $I$ has finitely many conjugates in $B$.
\end{itemize}

\noindent In particular, $(B,+)$ and $(B,\circ)$ are finitely presented.
\end{theo}
\begin{proof}
Let $S=\operatorname{Soc}(B)$; so $B/S$ is a finite skew brace. Of course, if $a\in B$, then $b\ast a=\lambda_b(a)-a=0$ and $[a,b]_+=0$ for every $b\in S$. Let $n$ be the index of a normal abelian subgroup~$A$ of $(B,\circ)$. Let $C$ be the subset of $S$ made by all elements of the form~$nb$. Then $C$ is a characteristic subgroup of the additive group of $S$. On the other hand, if we look at $S$ as a subgroup of $(B,+)$, then we see that it is normal in the semi-direct product $G=(B,+)\rtimes_\lambda(B,\circ)$. This means that also $C$ is normal in $G$, whenever we see it as a subgroup of $(B,+)$. In particular, $C$ is $\lambda$-invariant and is a subgroup of $(B,\circ)$.

Note also that $C$ coincides with the set of all elements of $S$ of the form $b^n$, since $b^n=nb$, and consequently $C$ is contained in $A$. Moreover, Lemma 1.10 of \cite{Cedo} yields that \[
\begin{array}{c}\label{cedoeq}
\lambda_a(c)=a\circ c\circ a^{-1}\tag{$\dagger$}
\end{array}
\] for all $a\in B$ and $c\in S$. This means that $C$ is normal in $(B,\circ)$ and so is actually an ideal of $B$.

Now, $B/C$ is finite and so Lemma \ref{fgideal} yields that $C$ is finitely generated as an ideal of~$B$. Let $c_1,\ldots,c_\ell$ be elements of $B$ generating $C$ as an ideal. Note that every element of $A$ has finitely many conjugates in $(B,\circ)$, so \eqref{cedoeq} yields that $(c_1,0)\ldots,(c_\ell,0)$ have finitely many conjugates in $G$ and consequently that the normal subgroup $N$ generated by $(c_1,0),\ldots,(c_\ell,0)$ in $G$ is finitely generated; let $(d_1,0),\ldots,(d_m,0)$ be generators of this subgroup. 

Let $D=\{d\,:\, (d,0)\in N\}$. Then $D$ is a finitely generated central subgroup of $(B,+)$ (it is in fact additively generated by $d_1,\ldots,d_m$). Since $D$ is contained in $\operatorname{Soc}(B)$, it is also multiplicatively generated by $d_1,\ldots,d_m$. Since $D$ is $\lambda$-invariant, Equation \eqref{cedoeq} yields that $D$ is normal in $(B,\circ)$ and is therefore an ideal. Thus $D=C$.

It now follows from Theorem \ref{thfp} that $B$ is finitely presented. Finally, for every $a\in A$ and $b\in B$, we have  $b\ast a=\lambda_b(a)-a=b\circ a\circ b^{-1}-a$, which by \eqref{cedoeq} means that there are finitely many conjugates of type $b\ast a$.  Since $B/C$ is finite, the statement is proved.~\end{proof}

\medskip

In connection with Theorem \ref{teoooo}, Corollary \ref{corollario} and Theorem \ref{alternative}, we note that it is also possible to define the {\it upper socle series} replacing the  operator ‘‘$\operatorname{Ann}$’’ with the operator ‘‘$\operatorname{Soc}$’’ in the definition before Theorem \ref{teoooo}. This series of ideals is relevant in some nilpotency theory of skew braces: it turns out, for instance, that a skew brace~$B$ of nilpotent type is right nilpotent if and only if $B=\operatorname{Soc}_n(B)$ for some positive integer~$n$ (see \cite[Lemmas 2.15 and 2.16]{Cedo}). Since the hypothesis in each of the three quoted theorems can be written as $B=\overline{\operatorname{Ann}}(B)$, it is reasonable to ask if in some situations we can replace this hypothesis with the new one $B=\overline{\operatorname{Soc}}(B)$. It turns out that under the additional assumption that every element of $(B,\circ)$ has finitely many conjugates in $(B,\circ)$, the statements corresponding to Theorem \ref{teoooo}, Corollary \ref{corollario} and~The\-o\-rem~\ref{alternative} hold: it is enough to repeat their proofs using the argument given in~The\-o\-rem~\ref{nilpotent2}.

\section{Declarations}

\noindent{\bf Conflict of interest}\quad The authors declare that they have no conflict of interest.

\smallskip

\noindent{\bf Data availability}\quad Data sharing is not applicable to this article as no new data were created or analyzed in this study.

\bigskip\bigskip\bigskip

\renewcommand{\bibsection}{\begin{flushright}\Large
{
REFERENCES}\\
\rule{8cm}{0.4pt}\\[0.8cm]
\end{flushright}}

\bigskip\bigskip


\begin{flushleft}
\rule{8cm}{0.4pt}\\
\end{flushleft}

{
\sloppy
\noindent
Marco Trombetti

\noindent 
Dipartimento di Matematica e Applicazioni ``Renato Caccioppoli''

\noindent
Università degli Studi di Napoli Federico II

\noindent
Complesso Universitario Monte S. Angelo

\noindent
Via Cintia, Napoli (Italy)

\noindent
e-mail: marco.trombetti@unina.it 

}


\end{document}